\documentclass[11pt]{amsart}
\hoffset= - 0.75 in

\usepackage{amsfonts, amsmath, amssymb, amsthm}
\usepackage{latexsym}
\newtheorem{thm}{Theorem}

\newtheorem{lem}{Lemma}

\newcommand{\RR}{\mathbb{R}}

\newcommand{\bin}[2]{
	\left(
		\begin{array}{@{}c@{}}
			#1  \\  #2
		\end{array}
	\right)		}

\begin{document}

\title{On LP-orientations of cubes and crosspolytopes}
\author{Mike Develin}
\address{Department of Mathematics, UC-Berkeley, Berkeley, CA 94720-3840}
\curraddr{2706B Martin Luther King Jr. Way, Berkeley, CA 94703}
\date{\today}
\email{develin@math.berkeley.edu}

\begin{abstract}

In a paper presented at a 1996 conference, Holt and Klee introduced a set of necessary conditions for an orientation
of the graph of a $d$-polytope to be induced by a realization into $\RR^n$ and linear functional on that space. In
general, it is an open question to decide whether for a polytope $P$ every orientation of its graph satisfying these
conditions can in fact be realized in this fashion; two natural families of polytopes to consider are cubes and
crosspolytopes. For cubes, we show that, as $n$ grows, the percentage of $n$-cube Holt-Klee orientations which can be
realized goes asymptotically to 0. For crosspolytopes, we give a stronger set of conditions which are both necessary
and sufficient for an orientation to be in this class; as a corollary, we prove that all shellings of cubes are
line shellings.

\end{abstract}

\maketitle

\section{Introduction}

Given a polytope in affine space and a generic linear function $f$, we may define an orientation on the
graph of the polytope by orienting an edge $vw$ from $v$ to $w$ if and only if $f(v) < f(w)$. If we fix
the polytope $P\subset \RR^d$, a natural issue to investigate is which orientations of its graph are
induced by linear functions. This problem falls under the umbrella of linear programming.

From a combinatorial perspective, though, this question is too narrow. Two polytopes are defined to be 
{\em combinatorially equivalent} if they have the same face lattice. For a combinatorialist, a more
natural and broader question is the following: given a polytope $P$ and an orientation $O$ on its edge
graph $G(P)$, does there exist any combinatorially equivalent polytope and linear function $f$ which
defines this orientation? The realization space of combinatorially equivalent polytopes may be quite
complicated, and in general this is a very difficult question.

From now on, whenever we refer to a polytope $P$ we are referring to the equivalence class of all combinatorially
equivalent polytopes; by a {\em realization} we mean any polytope $Q\subset \RR^d$ that is in this equivalence class.
Following~\cite{HK}, we say that an orientation $O$ of the graph of a polytope $P$ is an {\em LP-orientation} if there
exists a realization of $P$ and a linear functional $f$ such that the orientation defined by $f$ on the graph of $P$
is precisely $O$. Directed graphs which can be expressed in this fashion are called {\em polytopal digraphs}; since we
are fixing the graph and the polytope, we prefer a term which characterizes just the orientation. Since many polytopes
may have the same graph, an orientation of the graph of $P$ may be a polytopal digraph (realized by some other
polytope) without being an LP-orientation with respect to $P$; whenever we use the term LP-orientation, we implicitly
fix the polytope $P$.

In \cite{HK}, Holt and Klee presented three necessary conditions for an orientation of $P$ (by which we
always mean an orientation of the graph of $P$) to be an LP-orientation. First of all, the induced
orientation on every face of the polytope must have a unique sink and unique source. Second of all, the
orientation must be acyclic; there cannot be any directed cycles. Finally (and this is a very
nontrivial condition), for each $k$-face of the polytope, there must be at least $k$ monotone ascending
paths from the relative source to the relative sink of the face, pairwise disjoint except at the
endpoints. We say that an orientation of a polytope is {\em Holt-Klee} if it satisfies these
conditions.

For some polytopes, these conditions are sufficient to check whether an
orientation is linearly inducible or not. For instance, Mihalisin and Klee
\cite{MK} showed that for 3-polytopes, all Holt-Klee orientations are
LP-orientations, and the conditions clearly suffice for simplices.  Two
natural families of polytopes to consider are cubes and crosspolytopes;
the result of Mihalisin and Klee shows that dimension 4 is the threshhold
for counterexamples.  Walter Morris~\cite{WM} has recently found a
Holt-Klee orientation of the 4-cube which is not a LP-orientation; in the
course of this paper, we will construct a counterexample for the
4-crosspolytope.

For cubes, we show that the Holt-Klee conditions are very insufficient. To be precise, as $n$ gets
large, the proportion of Holt-Klee orientations that are LP-orientations dwindles rapidly to 0. In
particular, the number of LP-orientations grows at most exponentially in $n$, while the number of
Holt-Klee orientations grows at least doubly exponentially in $n$. Our method is nonconstructive, and
the ideas used in the proof are general ones which we believe can be applied to investigate the
question for a wide variety of polytopes.

For crosspolytopes, we give a set of necessary and sufficient conditions for an orientation of a
crosspolytope to be an LP-orientation. These conditions are stronger than the Holt-Klee conditions, but
most Holt-Klee orientations of crosspolytopes are LP-orientations. The proof is entirely different from
the case of cubes, and relies heavily on the specific combinatorial structure of crosspolytopes.
 
\section{Cubes}

For cubes, we show that as $n$ gets large, almost no Holt-Klee orientations of the hypercube are
linearly inducible. Our method will be indirect; we will obtain an upper bound on the number of
linearly inducible orientations, and then construct a class of Holt-Klee orientations whose size grows
far faster than this upper bound.

\subsection{The number of LP-orientations of the $n$-cube}\label{linind}

In this section, we show that the number of LP-orientations of the $n$-cube grows at most
exponentially in $n$. We assume without loss of generality that the linear functional in
question is just $x_1$, and consider all realizations of the $n$-cube into $\RR^n$ and the
resulting orientations. To enumerate these possible orientations, we will use the following
result, a special case of a theorem by Basu, Pollack, and Roy~\cite{BPR}.

\begin{thm}[Basu, Pollack, Roy]\label{signcond}
Suppose we have a finite subset of polynomials $\{p_1,\ldots,p_s\}\subset
\RR[x_1,\ldots,x_k]$, each of degree at most $d$. For each point $x\in
\RR^k$, define the function $\sigma(x)$ to be the tuple of sign vectors
$(\text{sign }p_1(x), \ldots, \text{sign }p_s(x))$. Then the number of
distinct values of $\sigma(x)$ is at most $\bin{s}{k} O(d)^k$. 
\end{thm}

We first make the following key observation: a realization of the $n$-cube into $\RR^n$ is completely 
determined by the $2n$ facet-defining inequalities $f_i(x)\le w_i$, where we have $f_i\in 
(\RR^n)^\star$, $w_i\in \RR$.

It is certainly not the case that every system of facet-defining inequalities forms a
combinatorial $n$-cube, but this observation allows us to parametrize the set of all
combinatorial cubes as a subset of such $2n$-tuples of inequalities. Each such inequality has
$n+1$ real coefficients, so the set of realizations of the $n$-cube into $\RR^n$ (picking an
arbitrary order for the facet-defining inequalities) can be parametrized by a subset of
$\RR^{2n(n+1)}$.

Next, we wish to deduce from these $2n(n+1)$ parameters the orientation of each edge of the cube.
Consider any vertex $v$ of the cube. This is the unique point of intersection of the $n$ facets it is
contained in; it is the unique solution point of a matrix equation $Av=w$, where each row of the
equation states that the facet-defining inequality $f_i(x)\le w_i$ achieves equality at the vertex $v$,
and the $n$ rows correspond to the $n$ facets containing $v$. By Cramer's rule, we can express the
first coordinate of $v$ as the ratio of two determinants, which are polynomials of degree $n$ in our
$2n(n+1)$ parameters, say $g_v(t)/h_v(t)$, where $t\in \RR^{2n(n+1)}$.

The crucial point here is that if the facet-defining inequalities form a cube, then this point
$v$ exists and is unique. Thus, in the subset of $\RR^{2n(n+1)}$ we are considering, the
denominator determinant will be nonzero. The orientation of an edge between two vertices $v$
and $w$ is determined by the sign of the expression $g_v(t)/h_v(t) - g_w(t)/h_w(t)$. This
expression is well-defined if the facet-defining inequalities form a cube; furthermore, it is
nonzero if the cube in question actually has an orientation induced by the linear functional
$x_1$ (i.e. if no two adjacent vertices have equal $x_1$ coordinates.)

This sign is the product of the signs of $g_v(t)h_w(t)-g_w(t)h_v(t)$ and $h_v(t)h_w(t)$; both
of these signs are nonzero in the subset of $\RR^{2n(n+1)}$ which corresponds to combinatorial
cubes by previous remarks, but also both of these polynomials have degree at most $2n$ and are
defined over all of $\RR^{2n(n+1)}$.

So the orientation of a combinatorial cube is completely determined by the signs of these polynomials
as $(v,w)$ ranges over all pairs of adjacent vertices of the cube. We now apply Theorem~\ref{signcond}
to bound the number of sign conditions. There are $n(2^n)$ edges in the $n$-cube, so there are
$2n(2^n)$ polynomials. Each has degree at most $2n$, and the dimension of the space in question is
$2n(n+1)$. Since the number of sign conditions on a subset of $\RR^{2n(n+1)}$ is at most the number of
sign conditions on the entire space, we obtain that the number of LP-orientations of the $n$-cube is at
most $\bin{2n(2^n)}{2n(n+1)} O(2n)^{2n(n+1)}$. Since $\bin{n}{k} \le n^k$, some straightforward if
messy algebra shows that this number is bounded by an expression of the form $p(n)^{q(n)}$, where
$p$ and $q$ are polynomials in $n$.

\subsection{A large class of Holt-Klee orientations for the $n$-cube}\label{hk}

In this section, we show that the number of Holt-Klee orientations of the $n$-cube grows at
least doubly exponentially in $n$. We accomplish this by exhibiting a large class of Holt-Klee
orientations of the $n$-cube which grows this fast with respect to $n$.

Given $n$, we define this class of orientations as follows. The edges of the $n$-cube (which is
assumed to have vertices labeled with binary strings of length $n$) can be canonically
partitioned into $n$ sets $\{E_1,\ldots,E_n\}$ as follows: the two vertices at the endpoint of
an edge differ in one of the $n$ places, and if it is place $i$, put the edge into set $E_i$.

Given an edge $e$, we orient it as follows. If $e$ is in $E_i$ for $i\neq n$, we orient it
towards its endpoint with 1 in place $i$. If $e$ is in $E_n$, let $r=\lfloor
\frac{n}{2}\rfloor$, and consider the substring consisting of the first $n-1$ places, which is
the same for both of the endpoints of $e$. If this string has fewer than $r$ 1's, orient this
edge towards the endpoint with 1 in the last place. If it has greater than $r$ 1's, orient it
towards the endpoint with 0 in the last place.

We have now filled in nearly the entire orientation of the $n$-cube. We will show that we can fill
in the remaining $\bin{n-1}{r}$ edges arbitrarily, and that the resulting 
orientation of the
$n$-cube will be Holt-Klee.

Take any such orientation. The orientation is clearly acyclic, since every edge not in $E_n$ is
oriented towards the endpoint with more 1's. Each face has a unique sink (it must be in the set
of vertices with the largest number of 1's among the first $n-1$ coordinates, which is either 1
or 2, and if it's 2 the orientation of the $E_n$-edge between them gives it to us) and unique
source. The only nontrivial thing to check is the path condition: in each $k$-face, we must
have $k$ paths between source and sink, disjoint except at the endpoints.

By the {\em standard k-cube}, we mean the $k$-cube which has all edges oriented towards the
endpoint with an additional 1; this is easily shown to satisfy the path condition By design,
each $k$-face not including edges from $E_n$ is isomorphic to the standard $k$-cube and thus
satisfies the path condition.

In faces $F$ which do contain edges from $E_n$ (and thus one incident on every vertex), as
always, all of the edges not from $E_n$ are oriented towards the vertex with an additional 1.
By the way we've oriented edges in $E_n$, one of two cases must hold.

\textit{Case 1.} The edges in $F$ belonging to $E_n$ are all oriented in the same direction. In
this case, $F$ is isomorphic to the standard $k$-cube, and thus satisfies the path condition.

\textit{Case 2.} The source and sink are both vertices with final coordinate 0. In this case,
we can find $k-1$ disjoint (except at the endpoints) paths contained inside the subcube of
dimension $k-1$ consisting of all vertices with final coordinate 0, since this is just
isomorphic to the standard $(k-1)$-cube (containing no edges from $E_n$.) 
For the final path, we
simply start with the edge in $E_n$ incident on the source, then find a path in the
$(k-1)$-subcube consisting of all vertices with final coordinate 1 (also 
isomorphic to the
standard $(k-1)$-cube) from its relative source (which has the same first 
$n-1$ coordinates as
the source of the $k$-cube) to its relative sink (which has the same first $n-1$ coordinates as
the sink of the $k$-cube), and end with the edge in $E_n$ incident on the sink. This produces
$k$ paths from source to sink, which are disjoint except at the endpoints. So $F$ satisfies the
path condition in this case as well.

Therefore, we have shown that any orientation in our class is in fact a Holt-Klee orientation. The size
of this class is $2^{\bin{n-1}{r}}$. $\bin{n}{n/2}$ grows as $2^n/\sqrt{n}$, so the size of this class
grows as $2^{2^n/\sqrt{n}}$, which is doubly exponential in $n$. In particular, the number of Holt-Klee
orientations is at least $2^{2^n/r(n)}$, where $r$ is a polynomial in $n$.

\subsection{Main Result}
The main result now follows easily.

\begin{thm}
For large enough $n$, there exist Holt-Klee orientations of the $n$-cube which are not linearly
inducible. Indeed, as $n$ grows, the fraction of Holt-Klee orientations of the $n$-cube which
are linearly inducible approaches 0 rapidly.
\end{thm}

\begin{proof}

This follows from the two previous sections. The number of Holt-Klee
orientations of the $n$-cube is at least $2^{2^n/r(n)}$, while the number
of linearly inducible orientations is at most $p(n)^{q(n)}$ ($p(n), q(n),
r(n)$ are all polynomials in $n$.) We claim that for large enough $n$,
$2^{2^n/r(n)} > p(n)^{q(n)}$ for any polynomials $p,q,r$; indeed, taking
logs base 2 of both sides yields that this is equivalent to $2^n/r(n) >
q(n)\text{log }p(n)$ or $2^n > r(n) q(n) \text{log }p(n)$. But the 
right
side is (bounded by) a polynomial in $n$, while the left side is
exponential in $n$, so for sufficiently large $n$ this is true, and in
fact the disparity grows quite quickly.

\end{proof}

The argument in Section~\ref{linind} provides a reasonable bound on the
number of linearly inducible orientations of any polytope; if a $d$-polytope has $e$ edges and
$f$ facets, the number of linearly inducible orientations is at most $\bin{2e}{f(d+1)}
f^{f(d+1)}$. If we have a family of polytopes for which the number of facets is linear in the
dimension and the number of vertices is at most exponential, then this is bounded by
$p(n)^{q(n)}$ for polynomials $p$ and $q$ ($n$ being the dimension), so this approach should work well 
on the few-facet case. For the opposite end of the spectrum, when the polytope has few vertices, we 
refer the reader to Mihalisin ~\cite{Mpre}.

The second part of the technique, computing an appropriate class of Holt-Klee orientations with
sufficient size, is more specific to the individual problem at hand. For cubes, we were able to
exhibit a class whose size grows doubly exponentially in $n$; obviously, a necessary
prerequisite for this is that the number of edges of a family grows exponentially in $n$.
Enumerating Holt-Klee orientations is, in general, a difficult problem, but in the case of
polytopes which, like cubes, have a lot of structure, it is feasible to construct a large class
of orientations. To take one example, this proof can be easily modified to show that for any
polytope $P$, the Cartesian product polytope $P\times I^n$ will, for large enough $n$, also
have the property that almost no Holt-Klee orientations are linearly inducible.

\section{Crosspolytopes}
We first state the main theorem, which we will prove by developing a technique known as pair
encoding and attacking the problem in its new formulation.

\begin{thm}\label{main}

Let $O$ be an acyclic orientation of the $d$-dimensional crosspolytope,
whose vertices naturally come in pairs $P_1,\ldots,P_d$. Then $O$ is
linearly inducible if and only if no proper subset of the pairs forms an
initial set for the orientation; that is, for no subset $S$ of the 
vertices consisting of a union of pairs is it the case that all edges 
between vertices $x\in S$ and $y\notin S$ are oriented towards $y$.

\end{thm}

The Holt-Klee conditions for the crosspolytope are relatively trivial; all proper faces are simplicial,
so if we have an acyclic orientation we only need to check globally that it has a unique source and
sink. The path condition then follows rather straightforwardly.

An acyclic orientation of the crosspolytope can have multiple sources (or sinks) if and only if there
are exactly two, and they are in the same pair. This is tantamount to the existence of a subset $S$ as
above consisting of either 1 or $d-1$ pairs; the conditions in the theorem are therefore a
strengthening of the Holt-Klee conditions. In particular, the Holt-Klee conditions are insufficient to
ensure that an orientation is an LP-orientation.

\subsection{Pair encoding}

In this section, we develop notation for the problem at hand. Given an acyclic orientation of a
$d$-crosspolytope, we translate this bijectively into a partition of $\{1,\ldots,2d\}$ as follows.

As stated before, the vertices of the $d$-crosspolytope naturally come in pairs $P_1,\ldots,P_d$; the
proper faces of this simplical polytope are then simply the sets of vertices with at most one element
from each pair.  Given any acyclic orientation of any graph with $n$ vertices, we can find a labeling
of the vertices of the graph with ${1,\ldots,n}$, not necessarily unique, such that all edges are
oriented towards the vertex with the larger label. Now, given an acyclic orientation $O$ of a
$d$-crosspolytope, we define the {\em pair sequence} of $O$ to be the partition of $\{1,\dots,2d\}$
into $d$ pairs given by assigning a labeling of the vertices corresponding to $O$ and then considering
the pairs of labels given to the pairs of vertices $P_i$. We will always list the pairs in
ascending order of smaller element, for instance as $(14)(25)(36)$, which corresponds to an
orientation of the octahedron.

Given any acyclic orientation $O$, this pair sequence is unique. This is because the labeling is unique
up to switching two adjacent labels on vertices in the same pair, since every other pair of vertices
has an edge between them; this switch does not affect the pair sequence. Consequently, when determining whether or 
not an acyclic orientation is an LP-orientation, it suffices to consider the related pair sequence.

Given a pair in a pair sequence of length $d$, we may {\em eliminate} this pair to get a pair sequence
of length $d-1$. The methodology here is straightforward; the remaining numbers give an
ordering of $2d-2$ elements in $d-1$ pairs, and we simply renumber so that these labels consist of
$\{1,\ldots,2d-2\}$. For instance, eliminating the pair $(25)$ from $(14)(25)(36)$ initially
yields the sequence $(14)(36)$, which reduces to $(13)(24)$. It is clear that if the ordering
of the initial pairs follows our listing convention (ascending order of smallest element), then the
ordering of the resulting smaller pair sequence will as well.

\subsection{Proof of the Main Theorem}

With pair encoding, Theorem~\ref{main} can be rephrased as follows.

\begin{thm}\label{pairmain}

Given an acyclic orientation $O$ of a $d$-dimensional crosspolytope, the orientation is an
LP-orientation if and only if, for all $0<k<d$, the first $k$ pairs in the corresponding pair sequence
do not comprise the elements $\{1,\ldots,2k\}$.

\end{thm}

Before we proceed to the proof of the theorem, we state and prove a useful lemma.
\begin{lem}\label{induct}
Suppose that the vertices $v_1,\ldots,v_{2d-2}$ form a $d-1$-dimensional crosspolytope $P$ in $d$-dimensional space. Then for any $y,z$ not 
in the affine hull of $P$ for which the line segment $yz$ meets the relative interior of $P$, the convex hull of the vertices 
$v_1,\ldots,v_{2d-2}, y, 
z$ is a $d$-dimensional crosspolytope.
\end{lem}

\begin{proof}

The proof is fairly simple. Assume without loss of generality that the vertices $v_1,\ldots,v_{2d-2}$ lie in the
hyperplane $x_d=0$; by making an appropriate affine transformation, we can assume that $y_i=z_i$ for all coordinates
aside from the last one, where $y_d$ is positive and $z_d$ is negative. We denote the point of intersection of the 
line segment $yz$ with
the polytope $P$ by $w$; this point has the same first $d-1$ coordinates as $y$ and $z$ and has $w_d=0$. (Affine
transformations, of course, do not change the combinatorial class of the convex hull.)  Now, we
must show that the proper faces of the polytope $Q$ which is the convex hull of all $2d$ vertices are just $\{F\cup
Z\}$, where $F$ is a proper face of $P$ and $Z$ is either empty, $\{y\}$, or $\{z\}$.

First of all, we show that all of these are actually faces of $Q$. Take any proper face $F$ of $P$, and consider a
face-defining linear functional $f$ on $\RR^{d-1}$ which is 0 on the convex hull of the vertices of $F$ and positive
elsewhere; in particular, we have $f(w)>0$, since $w$ does not lie on any proper face of $P$. Then $f+0x_d$ defines
the face $F$ of $Q$, while $f-(f(w)/y_d)x_d$ defines the face $F\cup \{y\}$ and $f+(f(w)/z_d)x_d$ defines the face 
$F\cup 
\{z\}$.
So the faces of $Q$ certainly include this set.

Conversely, suppose we have a face $G$ of $Q$ which is not in this set. Let $f$ be a linear functional defining it, so that $f=0$ on $G$ and
$f>0$ on $Q-G$. Then restricting $f$ to the hyperplane $x_d=0$, $f$ must define a face of $P$; $P$ is a subset of $Q$, so $f$ must be
nonnegative no $P$. The face of $Q$ defined by $f$ will be the union of this face with the subset of $\{y,z\}$ on which $f=0$; since $G$ is not 
in the proposed set of faces of $Q$, this subset must be $\{y,z\}$. But then $f$ must also be 0 on the point $w$, which implies (since it is 
nonnegative and $w$ is in the relative interior of $P$) that it is 0 on all of $P$. So $f$ must be 0 on all of $Q$, and hence $G=Q$ is not a 
proper face.
\end{proof}

\begin{proof}[Proof of Theorem~\ref{pairmain}]
For ease of notation, we will call a pair sequence which satisfies the condition of the theorem a {\em 
good} pair sequence, and one that  does not a {\em bad} pair sequence.

We will prove each direction of this theorem independently: first, we will show that if we have a good pair sequence, then the corresponding
orientation is an LP-orientation. The proof is by induction on $d$. For $d=1$ the claim is clear, since the only pair sequence is $\{12\}$, 
which is good and clearly an LP-orientation.

Given a pair sequence of length $d>1$, we simply eliminate the last pair; we claim that the resulting pair sequence
will also be good. Suppose not; then the first $k$ pairs of the smaller sequence, for some $1<k<d-1$, consist of the
numbers $\{1,\ldots,2k\}$. However, this means that there is some pair $P$ which has none of these numbers, and hence
has smallest element greater than any element in any of these pairs. Since we removed the last pair, the smallest
element of the removed pair must be greater than the smallest element of $P$, and hence both elements of the removed
pair must be larger than $2k$ in the original ordering. Therefore, the first $k$ pairs of the smaller sequence are
also the first $k$ pairs of the original sequence, which is a contradiction since the original sequence was good.
Since this is impossible, the resulting pair sequence must be good as well.

The remaining step in this direction is to show that if we have an LP-orientation corresponding to the smaller
sequence, and the original sequence was good, then we also have an LP-orientation corresponding to the original
sequence. Let the pair we're trying to add be $(lm)$. Since the original sequence is good, this pair is not
$((2d-1)(2d))$.

Suppose without loss of generality that the LP-orientation corresponding to the smaller sequence is witnessed by the
polytope $P\subset \RR^{d-1}$ and the linear function $f = x_1$, and embed that space $\RR^{d-1}$ into $\RR^d$ via
$(x_1,\ldots,x_{d-1})\rightarrow (x_1,\ldots,x_{d-1}, 0)$. Then the linear function $x_1$ still induces the same
ordering on all the vertices except the removed pair; our job is to add the removed pair in such a way that the
resulting polytope is a $d$-crosspolytope, and such that the linear function $x_1$ inserts the two new vertices at
spots $l$ and $m$ in order to correspond to the original pair sequence.

For this, we use Lemma~\ref{induct}. First, we find specific distinct values for $l_1$ and $m_1$ for which adding two
points with those $x_1$ coordinates induces the desired pair sequence. Because $l$ and $m$ are neither the lowest two
vertices nor the highest two vertices in the pair decomposition, we can then find a value $l_1<r_1<m_1$ such that
there exists a point $r$ in $P$ with first coordinate $r_1$.

Now, simply draw any line segment through $r$ which is contained neither in the hyperplane $x_d=0$ nor the hyperplane
$x_1=r_1$. Points on this line take all values of $x_1$; simply take the points $l$ and $m$ to have first
coordinates $l_1$ and $m_1$ respectively. Then $l$ and $m$ lie on opposite sides of the hyperplane $x_d=0$, since
$l_1<r_1<m_1$, and so the line segment between $l$ and $m$ intersects $P$ (at $r$)  as desires. Consequently, by
Lemma~\ref{induct}, the convex hull of the original $2d-2$ vertices of $P$ (which is a $d-1$-dimensional
crosspolytope) and the added vertices $l$ and $m$ is in fact a $d$-dimensional crosspolytope. This $d$-crosspolytope 
is then a realization showing that the original pair sequence corresponds to an LP-orientation. This completes the first 
direction of the proof.

For the other direction, we need to show that any bad pair sequence cannot correspond to an LP-orientation. Take
$1<k<d$ such that the first $k$ pairs of the sequence in question comprise the set $\{1,\ldots,2k\}$, and suppose we
have a linear function $f$ on a realization $P$ which induces this bad pair sequence. By shifting $f$ by a constant,
we can assume that $f(v_i)$ is negative for $i\le 2k$ and positive for $i>2k$. We now indulge the reader in some
algebra by which we produce a face which, combinatorially speaking, simply cannot be a face of the crosspolytope.

The basic process is simple; we start with our function which is negative on the first $k$ pairs and positive on all
the other vertices, and then modify it by adding multiples of face-defining linear functionals. Given a functional
$g$, we define the sign-pair sequence of $g$ to be the sequence of pairs of signs of $g$ on the two vertices of a pair
in the pair sequence in question. In particular, the sign-pair sequence of $f$ is $(--)(--)\cdots(--)(++)\cdots(++)$;  
there are $k$ pairs of $(--)$ and $d-k$ pairs of $(++)$. Our goal is to find a linear function with the sign-pair
sequence $(--)(++)\cdots(++)$ (after appropriate reordering of the pairs.)

Now, because of the combinatorial structure of the crosspolytope, if we take a set consisting of at most one member of each set,
we obtain a face. Consequently, we can find a linear functional $h$ with the sign-pair sequence $(0+)\cdots (0+)(++)\cdots(++)$,
where the first $k$ sign-pairs are $(0+)$ and the last $d-k$ are $(++)$. We consider the functions $f+\epsilon h$ for increasing
$\epsilon>0$ until one of the signs is different from that of $f$; the sign-pair sequence of the resulting function, after 
appropriate reordering of pairs and entries within a pair, will be $(0 0/-) (0/- 0/-) \cdots (0/- 0/-) (++) \cdots (++)$, where 
the $0/-$ entries reflect the fact that multiple entries may become 0 for the same value of $\epsilon$. 

Now, we simply repeat this process, except that we make sure that entries which are already be 0 are in the face we
use at the next stage, so that once a pair obtains a 0 it will keep that 0 for the rest of the process. We stop after
we have a 0 in every face, so that the resulting linear functional will have pair-sequence (where $x$ just represents
a sign we don't know) $(0 x)(0 x)\cdots (0 x)(0 x)(++)\cdots(++)$; that is, it has one zero in each of the first $k$
pairs, and is positive on both members of each of the last $d-k$ pairs.

In the last pair to be zeroed, necessarily, the sign of the element which did not get zeroed out must be nonpositive.  
(If it were positive, since it started out negative, it would have hit 0 along the way and retained that 0 through the
rest of the process.) Consequently, at least one of the $x$ values is not $+$. Now, we consider the face formed by the
$k$ zeroes in the expression above, and add multiples of a linear functional that defines that face, which has
sign-pair sequence $(0+)\cdots(0+)(++)\cdots(++)$. As we do this, all $x$'s which start out negative will eventually
become 0 for an instant and then positive; we stop at the moment at which the last negative sign becomes $0$ (which
may be at the beginning if the necessary nonpositive sign or signs are all 0.)

At this point, our linear function has sign-pair sequence (after reordering of pairs to bring the critical pair as above to the 
front) $(00)(0 0/+)\cdots(0 0/+)(++)\cdots(++)$. Since all signs are nonnegative, it defines a face of the crosspolytope. This 
face is not the entire crosspolytope, since it has at least one positive sign. But it contains two members of the same pair. We 
have reached a contradiction, since the combinatorial definition of a crosspolytope yields that it has no proper faces with two 
members of the same pair.

Consequently, if we have a bad pair sequence, it cannot possibly correspond to an LP-orientation. This completes the proof of 
Theorem~\ref{pairmain}.
\end{proof}

\subsection{Corollaries}

Theorem~\ref{pairmain} completely solves the question of determining whether an orientation of a crosspolytope is an
LP-orientation; the conditions in question are stronger than the Holt-Klee conditions. In particular, the pair
sequence $(13)(24)(57)(68)$ (among others) corresponds to an orientation of the 4-dimensional crosspolytope which
satisfies the Holt-Klee conditions but which is not an LP-orientation.

The map between acyclic orientations of the $d$-crosspolytope (given combinatorially as $d$ pairs of vertices) to pair
sequences is precisely $2^d(d!)$ to 1, corresponding to switching the elements of a pair and permuting the pairs.
Consequently, we can answer the question of how many LP-orientations the crosspolytope has by enumerating the number
of good pair sequences.  This is a fairly simple question; the number of good pair sequences of length $d$, which we
denote by $a_d$, can easily be computed recursively, via the equation $a_d = d!! - \sum_{k=1}^{d-1} k!! a_{n-k}$; this
drops out easily from looking at the first ``break point'' of bad pair sequences.

A more surprising corollary is that Theorem~\ref{pairmain} is equivalent to the statement that all shellings of the
$n$-cube are line shellings. Looking at the dual space, an ordering of the vertices corresponds to an ordering of the
facets of the $n$-cube; a moment's thought will establish that an ordering of the facets is a shelling if and only if
the ordering of the vertices corresponds to a good pair sequence. On the other hand, an ordering induced by a linear
function $f$ in a polytope $P$ corresponds precisely to the order of the facets in the line shelling of the polar
polytope $P^\Delta$ in dual space associated to the line defined by the point $f$. Consequently, translating
Theorem~\ref{pairmain} into dual space yields the statement that an ordering of the facets of a cube is a shelling if
and only if it is a line shelling.  Indeed, the statement that line shellings are actually shellings provides an
alternate proof that sequences associated to LP-orientations must be good, although we have chosen to present a direct
proof in order to demonstrate the general method.

\section*{Acknowledgements} 

The author was supported by an NSF Graduate Fellowship. I would like to thank Jed Mihalisin for bringing the problem
to my attention, for discussions on the subject and for formulating the conjecture that became Theorem~\ref{pairmain},
as well as for giving generous and useful comments on an early draft of this paper. I would also like to thank Saugata
Basu for referring me to the result on partitioning $\RR^n$ by hypersurfaces, and Bernd Sturmfels for noting a
simplification of the main proof technique in the case of hypercubes.


\begin{thebibliography}{99} 

\bibitem{BPR}
S. Basu, R. Pollack, M.-F. Roy, ``On the number of cells defined by a family of polynomials on 
a variety'', \textit{Mathematika} \textbf{43} (1996), 120--126

\bibitem{HK}
F. Holt and V. Klee, ``A proof of the strict monotone 4-step conjecture,'' \textit{Contemp. 
Math.} \textbf{223} (1999), 201--216

\bibitem{Mpre}
J. Mihalisin, ``Affine orientations of polytopes with few vertices,'' preprint.

\bibitem{MK}
J. Mihalisin and V. Klee, ``Convex and linear orientations of polytopal graphs,'' 
\textit{Discrete and Computational Geometry} \textbf{24}(2000), no. 2-3, 421--435 

\bibitem{WM}
W. D. Morris, Jr., ``Distinguishing cube orientations arising from linear 
programs,'' preprint.
\end{thebibliography}
\end{document}